\documentclass[letterpaper,12pt]{article}
\usepackage[english]{babel}
\usepackage[utf8x]{inputenc}
\usepackage[T1]{fontenc}
\usepackage{url}
\usepackage{setspace}
\usepackage[title, titletoc]{appendix}
\usepackage{yfonts}

\usepackage[letterpaper,top=3cm,bottom=2cm,left=3cm,right=3cm,marginparwidth=1.75cm]{geometry}

\usepackage{amsmath}
\usepackage{amssymb} 
\usepackage{csquotes}
\usepackage[colorlinks=true, allcolors=blue]{hyperref}
\usepackage{hyperref}
\usepackage{setspace}
\usepackage{graphicx}
\usepackage{listings}
\usepackage{epigraph}
\usepackage{wasysym}

\numberwithin{equation}{section}

\title{Quantum Electrodynamics (QED) Renormalization is a logical paradox, Zeta Function Regularization is logically invalid, and both are mathematically invalid}

\author{Ayal I. Sharon \thanks{Patent Examiner, U.S. Patent and Trademark Office (USPTO). ORCID 0000-0002-5690-6181. This research received no funding, and was conducted in the author's off-duty time. The opinions expressed herein are solely the author's, and do not reflect the views of the USPTO, the U.S. Dept. of Commerce, or the U.S. Government. } 
}

\begin{document}
\maketitle
\begin{abstract}

Quantum Electrodynamics (QED) renormalizaion is a paradox. It uses the Euler-Mascheroni constant, which is defined by a conditionally convergent series. But Riemann's series theorem proves that any conditionally convergent series can be rearranged to be divergent. This contradiction (a series that is both convergent and divergent) is a paradox in "classical" logic, intuitionistic logic, and Zermelo-Fraenkel set theory, and also contradicts the commutative and associative properties of addition. Therefore QED is mathematically invalid. 

Zeta function regularization equates two definitions of the Zeta function at domain values where they contradict (where the Dirichlet series definition is divergent and Riemann's definition is convergent). Doing so either creates a paradox (if Riemann's definition is true), or is logically invalid (if Riemann's definition is false). We show that Riemann's definition is false, because the derivation of Riemann's definition includes a contradiction: the use of both the Hankel contour and Cauchy's integral theorem. Also, a third definition of the Zeta function is proven to be false. The Zeta function has no zeros, so the Riemann hypothesis is a paradox, due to material implication and "vacuous subjects".

\end{abstract}

\pagebreak
\tableofcontents
\onehalfspacing

\section{Main Results}

\subsection{Quantum Electrodynamics (QED) Renormalization is a Paradox}

Richard Feynman called Quantum Electrodynamics (QED) renormalization a "dippy" process, and suspected that it is "not mathematically legitimate."
\footnote{See Feynman \cite{Feynman}, p.128: "The shell game that we play ... is technically called 'renormalization.' But no matter how clever the word, it is what I would call a dippy process! ... I suspect that renormalization is not mathematically legitimate."}
We show that QED renormalization contains a contradiction which renders it a paradox in "classical" logic, intuitionistic logic, and Zermelo-Fraenkel set theory, and renders it mathematically invalid. 

QED renormalizaion is a logical paradox because its use of the Gamma function includes use of the Euler-Mascheroni constant. This constant is defined by a conditionally convergent series. (More specifically, it is the difference between two divergent series). But Riemann's series theorem proves that the elements of any conditionally convergent series can be rearranged to result in a divergent series. 

In "classical" logic, intuitionistic logic, and Zermelo-Fraenkel set theory, this result (a series that is both convergent and divergent) is an impermissible paradox, because a statement cannot be simultaneously true and false (e.g. "Series X is convergent at domain value Y"). Moreover, in mathematics, this result contradicts the associative and commutative properties of addition. Therefore, QED renormalizaion is invalid logically and mathematically, because the Euler-Mascheroni "constant" is actually a paradox.

\subsection{Zeta Function Regularization is Either a Paradox or Invalid}

Moreover, Zeta function regularization contains a contradiction. It equates two different definitions of the Zeta function: the Dirichlet series definition, and Riemann's definition. It equates them at domain values where they contradict (where the former is divergent and the latter is convergent). The Dirichlet series definition is easily proven to be true. If also Riemann's definition is true, then Zeta function regularization is a paradox. Alternatively, if Riemann's definition is false, then Zeta function regularization is invalid.

We prove that Riemann's definition is false where it contradicts the Dirichlet series definition, because Riemann's definition is the result of a contradiction: use of both Hankel's contour and Cauchy's integral theorem. This contradiction invalidates Riemann's definition of Zeta, and invalidates every physics theory that assumes that Riemann's definition of the Zeta function is true, including Zeta function regularization, \footnote{See e.g. Hawking \cite{Hawking}, p.133, \S 1; and Matsui et al. \cite{Matsui}, Eq.7, and Eq.29.} and also
the Casimir effect, \footnote{See Dittrich \cite{Dittrich}, pp.30-34; Tong \cite{Tong}, pp.38-40; and Matsui et al. \cite{Matsui}, Eq.8.}
Quantum Electrodynamics (QED), \footnote{See Dittrich  \cite{Dittrich}, p.34; and Bavarsad et al. \cite{Bavarsad}, Abstract, Eq.50, and Appendix A.}
Quantum Chromodynamics (QCD); \footnote{See Dittrich \cite{Dittrich2}; Dittrich \cite{Dittrich}, p.34; and Arnold et al. \cite{Arnold}, Abstract. But also see criticism by Dirac \cite{Dirac}, and Bilal \cite{Bilal}, p.4.} 
Yang-Mills theory, \footnote{See Witten \cite{witten1991};  and Aguilera-Damia  \cite{Aguilera-Damia}.}
Supersymmetry (SUSY), \footnote{See Elizalde \cite{Elizalde}; and Bordag et al. \cite{Bordag}, \S1 and \S2.}
Quantum Field Theory (QFT), and \footnote{See Penrose  \cite{Penrose}, pp.656,678; Schnetz \cite{Schnetz}, \S E; and Cognola et al. \cite{Cognola}, \S 1.}
Bosonic String Theory. \footnote{See He \cite{He}; Veneziano \cite{Veneziano2}; Freund \cite{Freund}; Toppan \cite{Toppan}; Nunez \cite{Nunez}; pp.17-18; Tong \cite{Tong}, pp.39-40; and Bordag et al. \cite{Bordag}, \S2.} Also, a third definition of the Zeta function (that contradicts both the Dirichlet series definition and Riemann's definition) is false. And, because the Zeta function is exclusively defined by the Dirichlet series, it has no zeros. This renders the Riemann hypothesis a paradox, due to its "vacuous subjects" (the non-existent zeros) and material implication.

\section{QED Renormalization}

\subsection{The Riemann Series Theorem}

We begin with the definitions of \textit{absolutely} convergent series,  \textit{conditionally} convergent series, and divergent series:
\begin{itemize}
\item Infinite series $\sum a_n$ is divergent if $\sum a_n$ does not converge to a single value. Convergence and divergence are mutually exclusive characteristics. \footnote{See Hardy \cite{Hardy}, p.1.} 

\item Infinite series $\sum a_n$ is \textit{absolutely} convergent if $\sum a_n$ 
converges to a single value, and $\sum |a_n|$  converges to a single value. 

\item Infinite series $\sum a_n$ is \textit{conditionally} convergent if $\sum a_n$ 
converges to a single value, but $\sum |a_n|$ is divergent.
\end{itemize}

According to the Riemann series theorem (a.k.a. the Riemann rearrangement theorem): 

\begin{quotation}
By a suitable rearrangement of terms, a \textit{conditionally} convergent series may be made to converge to any desired value, or to diverge.
\footnote{See Weisstein \cite{Weisstein}, citing Bromwich \cite{Bromwich2}, p.74. See also Gardner \cite{Gardner2}, p.171; and Havil \cite{Havil}, p.102.} 
\end{quotation}
Here is one proof that a \textit{conditionally} convergent series can be rearranged to diverge:
\begin{quotation}
We can also rearrange the terms of any conditionally convergent series so that it will diverge. One such rearrangement is to pick positive terms to add to a million, then add on one negative term, then add on positive terms to reach a trillion, then add on another negative term, then add positive terms till we are beyond a googolplex, then add on a negative term \ldots
\footnote{See Galanor \cite{Galanor}. For more detailed proofs, see Bona et al. \cite{Bona}, Ch.9, Sec. 61, pp.120-121, Lemma 1 and Theorem 43.}
\end{quotation}

Therefore, according to the Riemann series theorem, a  conditionally convergent series is both convergent and divergent, depending upon the arrangement of its terms. This contradiction is a logical paradox in all logics whose axioms include the Law of Non-Contradiction (LNC) (e.g. "classical" logic, intuitionistic logic, and Zermelo-Fraenkel set theory), and thus is mathematically invalid, due to mathematics having the LNC as an axiom. In addition, in mathematics, this result contradicts the associative and commutative properties of addition, which provides another reason why such a series invalidates any mathematical theorem that uses it.

\subsection{The Gamma Function, Defined Using the Euler–Mascheroni Constant}

\subsubsection{One Definition of the Gamma Function}

One example of a \textit{conditionally} convergent series is the Euler–Mascheroni constant ($\gamma$). The Euler–Mascheroni constant appears in one definition of the Gamma function. In regards to the Gamma function, Kar \cite{Kar} states:
\footnote{See Kar \cite{Kar}, p.6. See also Hochstadt \cite{Hochstadt}, Chapter 3, "The Gamma Function".}
\begin{quotation}
Divergent sums and integrals occur in mathematics and physics. To avoid them, one has to implement unintuitive methods to deduce finite values for divergent quantities.
The Gamma function serves as a good toy model for that ... The Euler constant $\gamma$ which is associated with the Gamma function illustrates how the difference of two divergent quantities can still lead to a finite value.
\end{quotation}

However, Kar's \cite{Kar} description identifies a crucial logical flaw - without recognizing it as such. The "unintuitive" methods implemented "to deduce finite values for divergent quantities" are not merely "unintuitive". They are in fact \textit{contradictions}, and therefore logically and mathematically  invalid. Assigning a \textit{finite} value to an \textit{infinite (divergent)} value is a contradiction that violates the Law of Non-Contradiction (LNC). The LNC is the most important axiom in Aristotelian logic, and is also an axiom of the "classical logic" of Whitehead and Russell's \textit{Principia Mathematica}, Brouwer and Heyting's intuitionistic logic, and Zermelo-Fraenkel set theory.

\subsubsection{The Euler–Mascheroni "Constant" is a Conditionally Convergent Series}

Kar's \cite{Kar} definition of the Gamma function includes the Euler–Mascheroni constant as a critical element:
\footnote{See Kar \cite{Kar}, p.7.}
\begin{quotation}
By taking the logarithm of the [Gamma function] and by separating the divergent constant as $k$, we get
\begin{equation}
\log \frac{1}{\Gamma(z)} = \log z +
\sum_{n=1}^{\infty} \log \Big(e^{−\frac{z}{n}} \Big[1+ \frac{z}{n} \Big]\Big) + zk.    
\end{equation}
The constant k can be determined by comparing the derivative of the logarithm
of the Gamma function to leading orders. We find that,
\begin{equation}
\psi(z) = \frac{d}{dz} \log \Gamma (z) \sim − \frac{1}{z} − \gamma + O(z).
\end{equation}
Thus the arbitrary constant $k$ is the Euler constant $\gamma$.
\end{quotation}

Kar \cite{Kar} then defines the Euler–Mascheroni constant ($\gamma$), as a conditionally convergent series with a finite limit:
\footnote{See Kar \cite{Kar}, p.8, Eq. 1.1.10. See also Weisstein \cite{Weisstein3}.}
\begin{equation*} \label{eq1}
\begin{split}
\gamma & = \lim_{n \to \infty} \Big( \sum_{i=1}^{n} \frac{1}{i} − \int_{1}^{n} \frac{dx}{x} \Big) \\
 & = \lim_{n \to \infty} \Big( \frac{1}{n} +
\sum_{i=1}^{n-1} \Big[ \frac{1}{i} − \log \big(1+ \frac{1}{i}\big)\Big] \Big) \\
 & = \sum_{i=1}^{\infty} \Big[ \frac{1}{i} − \log \Big(1+ \frac{1}{i}\Big)\Big] \\
 & = 0.57721 \ldots \\
 & = -\Gamma ' (1) \\ 
\end{split}
\end{equation*}

Kar \cite{Kar} then states that "Even though we have a divergent sum $\displaystyle \sum\limits_{i=1}^{n} \frac{1}{i}$ and a divergent integral $\displaystyle \int\limits_{1}^{n} \frac{dx}{x}$, the difference - the Euler constant - is finite."
\footnote{See Kar \cite{Kar}, p.8, Eq. 1.1.10.}
However, Kar's \cite{Kar} definition of the Euler-Mascheroni constant can easily be rewritten as follows:
\begin{equation*} \label{eq2}
\begin{split}
\gamma & = \sum_{i=1}^{\infty} \Big[ \frac{1}{i} − \log \Big(1+ \frac{1}{i}\Big)\Big] \\
 & = \sum_{i=1}^{\infty}   \Big( \frac{1}{i} \Big) + \sum_{i=1}^{\infty} \Big( − \log \Big(1+ \frac{1}{i}\Big)\Big) \\
\end{split}
\end{equation*}

When rewritten in this manner, the Euler-Mascheroni constant is clearly a conditionally convergent series. The series is \textit{conditionally} convergent because the sum of all terms in the series, $\sum \Big[ (i^{-1}) − \log \Big(1+ i^{-1}\Big)\Big]$, is convergent. But the sum of the absolute values of all terms, $\sum  |i^{-1}| + \sum | − \log \Big(1+ i^{-1}\Big)|$, is divergent.

The Riemann series theorem holds that any conditionally convergent series (e.g. the Euler-Mascheroni constant) can be rearranged to have any finite value, and can also be rearranged to be divergent. So the same series is both convergent and divergent. This is a contradiction. In logical terms, it is a paradox. This result also contradicts the associative and commutative properties of addition. So a conditionally convergent series introduces a contradiction into any mathematical "proof" that uses it.

\subsubsection{Rearranging the Euler–Mascheroni Series to be Divergent}

As discussed in the previous section, the conditionally convergent series that defines the the Euler-Maseroni constant can be rearranged to diverge. Here is one example of how to do so:
\begin{quotation}
We can also rearrange the terms of any conditionally convergent series so that it will diverge. One such rearrangement is to pick positive terms to add to a million, then add on one negative term, then add on positive terms to reach a trillion, then add on another negative term, then add positive terms till we are beyond a googolplex, then add on a negative term \ldots
\footnote{See Galanor \cite{Galanor}.}
\end{quotation}

\subsection{QED Renormalization and the Euler-Mascheroni Constant}

Andrey \cite{Andrey} discloses an example use of the Euler-Mascheroni constant in "QED at One Loop":
\footnote{See Andrey \cite{Andrey}, p.23.}

\begin{quotation}
In order to define a dimensionless coupling $\alpha$, we have to introduce a parameter $\mu$ with the dimensionality of mass (called the renormalization scale):
\begin{equation} \label{eq:renorm}
\frac{\alpha (\mu)}{4 \pi} = \mu^{-2\epsilon} \cdot \frac{e^2}{(4 \pi)^{d/2}} \cdot e^{-\gamma \epsilon}    
\end{equation}
where $\gamma$ [in the term $e^{-\gamma \epsilon}$] is the Euler constant. In practise, this equation is more often used in the opposite direction:
\begin{equation}
\frac{e_0^2}{(4 \pi)^{d/2}}  = \mu^{2\epsilon} \cdot \frac{\alpha (\mu)}{4 \pi} \cdot Z_{\alpha}(\alpha (\mu)) \cdot e^{\gamma \epsilon}    
\end{equation}
We first calculate some physical quantity in terms of the bare charge $e_0$, and then re-express it via the renormalized $\alpha (\mu)$.
\end{quotation}

Given the above discussion about the Euler-Mascheroni constant being a paradox that is mathematically invalid due to contradiction, the "QED at One Loop" is also a paradox that is mathematically invalid due to contradiction, due to its use of the Euler-Mascheroni constant. Andrey \cite{Andrey} further discloses in a footnote:
\footnote{See Andrey \cite{Andrey}, p.23.}

\begin{quotation}
The first renormalization scheme in the framework of dimensional regularization was called MS (minimal subtractions); in this scheme
\begin{equation*}
\alpha (\mu) = \mu^{-2\epsilon} \cdot \frac{e^2}{(4 \pi)}   
\end{equation*}
It soon became clear that results for loop diagrams in this scheme look unnecessarily complicated, and the $\overline{MS}$ (modified minimal subtractions) scheme [in Eq. (\ref{eq:renorm})] was proposed. Some authors use slightly different definitions in the $\overline{MS}$, with $\Gamma (1 + \epsilon)$ or $1/\Gamma (1 - \epsilon)$ instead of $e^{-\gamma \epsilon}$ [in Eq. (\ref{eq:renorm})].
\end{quotation}

When $\Gamma (1 + \epsilon)$ or $1/\Gamma (1 - \epsilon)$ are used instead of $e^{-\gamma \epsilon}$), it is not clear which definition of the Gamma function is being used. According to Emil Artin \cite{Artin}, one definition of the Gamma function (derived by Weierstrass) uses the Euler-Mascheroni constant.
\footnote{See Artin \cite{Artin}, pp.15-16, including Eq. (2.8).} 

Moreover, another definition of the Gamma function (that Artin \cite{Artin} attributes to Gauss) does not include the Euler-Mascheroni constant, but \textit{does include} a logical error. Artin's Eq.(2.7), which is proved for the interval $0<x\le1$, is:
\begin{equation}
\Gamma(x) = \lim_{n \to \infty}\frac{n^x n!}{x(x+1) \cdots (x+n)}
\end{equation}
Artin states:\footnote{See Artin \cite{Artin}, p.15. } "As $n$ approaches infinity, if the limit in Eq.(2.7) exists for a number $x$, it also exists for $x+1$." 
\textbf{However, no limit exists, because the fraction $\infty/\infty$ is undefined!!!} 

More specifically, if $x$ is Real and \textit{non-negative}, the product $n^x n!$ approaches infinity as $n$ approaches infinity, as does the product $x(x+1) \cdots (x+n)$. So the above definition of the Gamma function results in the undefined ratio $\infty/\infty$. In the alternative, if $x$ is Real and \textit{negative}, the product $n^x n!$ approaches the undefined ratio $\infty/\infty$ as $n$ approaches infinity. Also the product $x(x+1) \cdots (x+n)$ approaches infinity as $n$ approaches infinity. So the above definition of the Gamma function results in the undefined nested ratios $(\infty/\infty)/ \infty$. 
 
\section{Zeta Function Regularization}

\subsection{Two Contradictory Definitions of the Zeta Function}

Bernhard Riemann's famous paper \textit{On the Number of Primes Less Than a Given Magnitude} begins with the statement 
\footnote{See Riemann \cite{riemann1859number}, p.1.}
that the Dirichlet series definition of the Zeta function is "invalid" for all values of complex variable $s$ in half-plane $\text{Re}(s)\le1$. 

However, Riemann's use of the word "invalid" is flat-out \textit{wrong} in the context of formal logic. The fact that the Dirichlet series definition of the Zeta function is "divergent" in said half-plane does \textit{not} mean that the definition is logically or mathematically false in that half-plane, or that its derivation is logically "invalid".

The proof that the Dirichlet series definition of the Zeta function is "divergent" in said half-plane is not provided in Riemann's paper, but is easily found elsewhere.
\footnote{See also Hardy \cite{Hardy2}, pp.3-5, citing Jensen \cite{Jensen}, Cahen \cite{Cahen}, and Bromwich \cite{Bromwich}.}
\footnote{See also Hildebrand \cite{Hildebrand}, pp.117-119, Thm 4.6.}  
In addition, the "Integral Test for convergence" (a.k.a. the Maclaurin–Cauchy test for convergence) 
\footnote{The "Integral Test for convergence" is often taught in introductory calculus textbooks, to prove that the famous "harmonic series" is divergent. See e.g. Guichard \cite{Guichard1}, Thm 13.3.4.} 
proves that the Dirichlet series of the Zeta function is divergent for all values of $s$ on the Real half-axis $(\text{Re}(s)\le1, \text{Im}(s)=0)$, which is a sub-set of the half-plane of divergence ($\text{Re}(s)\le1$).

Moreover, the Dirichlet series of the Zeta function is also proven to be divergent for all values of $s$ on the misleadingly-named "line of convergence" ($\text{Re}(s)=1$), which is a sub-set of the half-plane of divergence, and which is the border line between the half-plane of divergence and the half-plane of convergence. 
\footnote{See Hardy \cite{Hardy2}, p.5, Example (iii), citing Bromwich  \cite{Bromwich}.}  
At the point $s=1$, the Dirichlet series of the Zeta function is the famous "harmonic series", which is proven divergent by the "Integral test for divergence". 
\footnote{See Guichard  \cite{Guichard1}, Thm 13.3.4.}
At all other values of $s$ on the "line of convergence", the Dirichlet series of the Zeta function is a bounded oscillating function, which by definition is divergent. 
\footnote{See Hardy \cite{Hardy2}, p.5, Example (iii), citing Bromwich \cite{Bromwich}.}  

Later in Riemann's paper, his so-called "analytic continuation" of the Zeta function 
\footnote{Note that Riemann does not use the expression "analytic continuation". Also, note that Riemann's method is very different from Weierstrass's "unit disk" method of analytic continuation. See Solomentsev \cite{Solomentsev}.}
results in a second definition of the Zeta function, one that Riemann claimed "always remains valid" (except at the point $s=1$). \footnote{See Riemann \cite{riemann1859number}, p.1.}    
In other words, Riemann's definition of the Zeta function is convergent for all values in both half-planes (except at the point $s=1$).

However, Riemann's terminology again confuses the logical concepts of "validity" and "invalidity" with the mathematical concepts of "convergence" and "divergence".
Riemann intended to claim that his definition of the Zeta function is "convergent" for all values of $s$ in half-plane $\text{Re}(s)\le1$ (except at $s=1$). 

However, this claim raises the issue of logical "validity". If Riemann's claim is true, \textbf{and if Zeta function regularization is valid,} then all of the following propositions are true:

Pair 1:
\begin{itemize}
    \item The Zeta function is divergent for all $s$ in half-plane $\text{Re}(s)\le1$.
    \item The Zeta function is convergent for all $s$ in half-plane $\text{Re}(s)\le1$, (except at $s=1$).
\end{itemize}

Pair 2:
\begin{itemize}
    \item The Zeta function is divergent for all $s$ on the Real half-axis, $s<1$.
    \item The Zeta function is convergent for all $s$ on the Real half-axis, $s<1$.
\end{itemize}

Pair 3:
\begin{itemize}
    \item The Zeta function is divergent for all $s$ on the "line of convergence" $\text{Re}(s)=1$.
    \item The Zeta function is convergent for all $s$ on the "line of convergence" $\text{Re}(s)=1$, (except at $s=1$).
\end{itemize}

If all of the above propositions are true, then the two contradictory definitions of the Zeta function form a logical paradox. However, contrary to Riemann's characterization, the proof that the Dirichlet series of the Zeta function is divergent throughout half-plane $\text{Re}(s)\le1$, is a logically \textit{valid} proof. But the proof of his definition is \textit{not valid}.

The divergence of the Dirichlet series definition of the Zeta function throughout the half-plane $\text{Re}(s)\le1$ does \textit{not} render the proof "invalid", \textit{nor} does it render the function false (\textit{nor} does it render the function "not valid", as per Riemann's incorrect terminology).
\footnote{In the nomenclature of logic, "valid" and "invalid" apply to arguments. "True" and "false" apply to propositions. A mathematical proof is an argument, and a mathematical function is a proposition.} 
In fact, it is Riemann's definition of the Zeta function that introduced a problem of logical invalidity into mathematics,  because Riemann's definition of Zeta contradicts the  Dirichlet series definition (which is proven to be divergent throughout the half-plane $\text{Re}(s)\le1$).  

In Riemann's defense, his paper \textit{On the Number of Primes Less Than a Given Magnitude} (1859) predates Frege's \textit{Begriffsschrift} 
\footnote{See Frege \cite{Frege3}.}  
(1879) by two decades, and predates the subsequent developments in logic and the foundations of mathematics by at least a half-century. Brouwer's \textit{The Untrustworthiness of the Principles of Logic} (1908), Whitehead and Russell's \textit{Principia Mathematica} \footnote{See Whitehead \cite{Whitehead2}.} 
(1910), \L ukasiewicz's \textit{On Three-Valued Logic} (1920), and   Zermelo–Fraenkel set theory (1920's) were all published long after Riemann's untimely death (1866) at the age of 39. The only relevant publication in the field of logic that was contemporaneous with Riemann's work was Boole's \textit{The Laws of Thought} 
\footnote{See \cite{Boole}.} 
(1854), of which Riemann clearly was unaware. 

\subsection{Derivation of Riemann's Zeta Function} \label{derivation}

In the derivation of the Riemann Zeta function,  Riemann uses the following equation: 
\footnote{See Riemann \cite{riemann1859number}, p.1.} 
\begin{equation}
\int_{0}^{\infty}e^{-nx}x^{s-1}\,dx = \frac{\prod(s-1)}{n^s}
\end{equation}
On the left side of the equation, Riemann uses the equation
\footnote{See Edwards \cite{Edwards}, p.9, fn 1.} 
$\sum_{n=1}^{\infty} r^{-n} = (r-1)^{-1}$ to replace the term $e^{-nx}$ in the integral with the term $(e^{x}-1)^{-1}$. On the right side of the equation, Riemann introduces a summation (from $n = 1$ to $\infty$) for the term $1/n^{s}$, thereby obtaining:
\begin{equation}
\int_{0}^{\infty} \frac{x^{s-1}}{e^x-1}\,dx = \prod(s-1) \cdot\sum_{n=1}^{\infty} \frac{1}{n^s}
\end{equation}
The Dirichlet series definition of the Zeta function defines $\zeta(s) = \sum n^{-s}$, so the above equation is rewritten as:
\begin{equation} \label{Eq1}
\int_{0}^{\infty} \frac{x^{s-1}}{e^x-1}\,dx = \prod(s-1) \cdot \zeta(s)
\end{equation}
Next, Riemann considers the following integral:
\begin{equation}
\int_{+\infty}^{+\infty} \frac{(-x)^{s}}{(e^{x}-1)} \cdot \frac{dx}{x}
\end{equation}
Edwards \cite{Edwards} states: 
\footnote{See Edwards \cite{Edwards}, p.10.}  
\begin{quotation}
The limits of integration are intended to indicate a path of integration which begins at $+\infty$ , moves to the left down the positive Real axis, circles the origin once once in the positive (counterclockwise) direction, and returns up the positive Real axis to $+\infty$. The definition of $(-x)^s$ is $(-x)^s = \exp[s\cdot \log(-x)]$, where the definition of $\log(-s)$ conforms to the usual definition of $\log(z)$ for $z$ not on the negative Real axis as the branch which is Real for positive Real $z$; thus $(-x)^s$ is not defined on the positive Real axis and, strictly speaking, the path of integration must be taken to be slightly above the Real axis as it descends from $+\infty$ to $0$ and slightly below the Real axis as it goes from $0$ back to $+\infty$.
\end{quotation}
 
This is the Hankel contour. 
\footnote{See Edwards \cite{Edwards}, pp.10-11; See also Whittaker  \cite{Whittaker}, pp.85-87, 244-45 and 266.}  
The first use of this contour integral path was by Hankel, in his investigations of the Gamma function. 
\footnote{See Weisstein \cite{Hankel_Contour}, citing Krantz \cite{Krantz}, \S13.2.4, p.159; and Hankel \cite{Hankel}.}

When the Hankel contour is split into three terms, it is written mathematically as follows. 
\footnote{See Edwards \cite{Edwards}, p.10.} 
The first term is "slightly above" the Real axis as it descends from $+\infty$ to $\delta$, the middle term represents the circle with radius $\delta$ around the origin, and the third term is "slightly below" the Real axis as it goes from $\delta$ back
 to $+\infty$.
\begin{equation} \label{Hankel}
\int_{+\infty}^{\delta} \frac{(-x)^{s}}{(e^{x}-1)}\cdot \frac{dx}{x} + \int_{|z|=\delta} \frac{(-x)^{s}}{(e^{x}-1)}\cdot \frac{dx}{x} + \int_{\delta}^{+\infty} \frac{(-x)^{s}}{(e^{x}-1)}\cdot \frac{dx}{x}
\end{equation}
In regards to the middle term (the circle term), Edwards \cite{Edwards} states: 
\footnote{See Edwards \cite{Edwards}, p.10.}  

\begin{quotation}
[T]he middle term is $2\pi i$ times the average value of $(-x)^s\cdot (e^{x}-1)^{-1}$ on the circle $|x|=\delta$ [because on this circle $i \cdot d \theta = (dx/x)$]. Thus the middle term approaches zero as $\delta \to 0$ provided $s>1$ [because $x(e^{x}-1)^{-1}$ is nonsingular near $x=0$]. The other two terms can then be combined to give[:]
\end{quotation}
 
\begin{equation}
\int_{+\infty}^{+\infty} \frac{(-x)^{s}}{e^{x}-1} \cdot \frac{dx}{x} = \lim_{\delta \to 0} \Big[ \int_{+\infty}^{\delta}  \frac {\exp[s(\log x - i\pi)]}{(e^{x}-1)}\cdot \frac{dx}{x} + \int_{\delta}^{+\infty} \frac{\exp[s(\log x + i\pi)]}{(e^{x}-1)}\cdot \frac{dx}{x} \Big]
\end{equation}
resulting in 
\begin{equation}
\int_{+\infty}^{+\infty} \frac{(-x)^{s}}{e^{x}-1} \cdot \frac{dx}{x} = (e^{i\pi s} - e^{-i\pi s})\cdot \int_{0}^{\infty} \frac{x^{s-1}\,dx}{e^{x}-1}
\end{equation}
Given that $(e^{i\pi s} - e^{-i\pi s}) = 2i\sin(\pi s)$, this can be rewritten as:
\begin{equation}  
\int_{+\infty}^{+\infty} \frac{(-x)^{s}}{e^{x}-1} \cdot \frac{dx}{x} = 2i\sin(\pi s)\cdot \int_{0}^{\infty} \frac{x^{s-1}\,dx}{e^{x}-1}
\end{equation}
Rearranging the terms results in:
\begin{equation}  \label{Eq2}
\int_{0}^{\infty} \frac{x^{s-1}\,dx}{e^{x}-1} = \frac{1}{2i\sin(\pi s)} \cdot \int_{+\infty}^{+\infty} \frac{(-x)^{s}}{e^{x}-1} \cdot \frac{dx}{x}
\end{equation}
The left sides of Equations \ref{Eq1} and \ref{Eq2} are identical, so Riemann equates the right sides of Equations \ref{Eq1} and \ref{Eq2}, resulting in Equation \ref{Eq3}:
\begin{equation}  \label{Eq3}
\int_{+\infty}^{+\infty} \frac{(-x)^{s}}{e^{x}-1} \cdot \frac{dx}{x} = 2i\sin(\pi s)\cdot \prod(s-1) \cdot \zeta(s)
\end{equation}
Then, Riemann multiplies both sides of the equation by $\prod(-s)\cdot s/ 2\pi is$, resulting in
\begin{equation}  
\frac{\prod(-s)\cdot s}{2\pi is} \cdot \int_{+\infty}^{+\infty} \frac{(-x)^{s}}{e^{x}-1} \cdot \frac{dx}{x} =  \frac{\prod(-s)\cdot s}{2\pi is} \cdot 2i\sin(\pi s)\cdot \prod(s-1) \cdot \zeta(s)
\end{equation}
The $s$ terms on the left side cancel out, as do the $2i$ terms on the right side, so
\begin{equation} \label{Eq4}
\frac{\prod(-s)}{2\pi i} \cdot \int_{+\infty}^{+\infty} \frac{(-x)^{s}}{e^{x}-1} \cdot \frac{dx}{x} =  \frac{\prod(-s)\cdot \prod(s-1) \cdot s}{\pi s} \cdot \sin(\pi s)\cdot \zeta(s)
\end{equation}
Next, 
\footnote{See Edwards \cite{Edwards}, p.8, Eq.5; and pp.421-425.}  the identity of 
$\prod(s) = s\cdot \prod(s-1)$ is substituted into Eq. \ref{Eq4}, resulting in:
\begin{equation} \label{Eq5}
\frac{\prod(-s)}{2\pi i} \cdot \int_{+\infty}^{+\infty} \frac{(-x)^{s}}{e^{x}-1} \cdot \frac{dx}{x} =  \frac{\prod(-s)\cdot \prod(s)}{\pi s} \cdot \sin(\pi s)\cdot \zeta(s)
\end{equation}
Finally, the identity 
\footnote{See Edwards  \cite{Edwards}, p.8, Eq. 6.} 
$\sin(\pi s) = \pi s\cdot \Big[\prod(-s)\prod(s)\Big]^{-1}$ is substituted into the right side of Eq. \ref{Eq5}, resulting in
\begin{equation} 
\zeta(s) = \frac{\prod(-s)}{2\pi i} \cdot \int_{+\infty}^{+\infty} \frac{(-x)^{s}}{e^{x}-1} \cdot \frac{dx}{x}
\end{equation}
This is the Riemann Zeta Function. 
\footnote{See Edwards \cite{Edwards}, pp.10-11. Eq.3.} 

\subsection{Riemann's Zeta Function is False, Because Hankel's Contour Contradicts Cauchy's Integral Theorem} 

However, as a reminder, in regards to the three terms of the Hankel contour shown in Equation \ref{Hankel}:
\footnote{See Edwards \cite{Edwards}, pp.10-11; and Whittaker \cite{Whittaker}, p.244-6, \S 12.22, citing Hankel \cite{Hankel}, p.7.}  
\begin{equation} 
\label{Hankel2}
\begin{aligned}
\int_{+\infty}^{+\infty} \frac{(-x)^{s}}{(e^{x}-1)} \cdot \frac{dx}{x} = & \int_{+\infty}^{\delta} \frac{(-x)^{s}}{(e^{x}-1)}\cdot \frac{dx}{x} 
& + \int_{|z|=\delta} \frac{(-x)^{s}}{(e^{x}-1)}\cdot \frac{dx}{x} 
& + \int_{\delta}^{+\infty} \frac{(-x)^{s}}{(e^{x}-1)}\cdot \frac{dx}{x}
\end{aligned}
\end{equation}

Edwards \cite{Edwards} states: 
\footnote{See Edwards \cite{Edwards}, p.10.} 
\begin{quotation}
[T]hus $(-x)^s$ is not defined on the positive Real axis and, strictly speaking, the path of integration must be taken to be slightly above the Real axis as it descends from $+\infty$ to $0$ and slightly below the Real axis as it goes from $0$ back to $+\infty$.
\end{quotation}
 
Riemann copied this solution directly from Hankel's derivation of the Gamma function $\Gamma(s)$. 
\footnote{See Whittaker \cite{Whittaker}, pp.244-5,266.} 
Riemann uses the Hankel contour in the derivation of the Riemann Zeta function. But neither Hankel nor Riemann provide an answer to the question that Edwards's comment leads to: 
\begin{quotation}
What is the mathematical basis for Hankel's "trick" of \textit{equating} the branch cut of $f(x)=\log(-x)$ to the limit of the Hankel contour ("slightly above" and "slightly below" the branch cut)? 
\end{quotation}

As every mathematician knows, the logarithm of a non-positive Real number is \textit{undefined}. So, by definition, all points on the branch cut have no defined value. Equating the branch cut to the limit of the Hankel contour ("slightly above" and "slightly below" the branch cut) is a \textit{de facto} assignment of values to points that, by the definition of logarithms, \textit{must have no value}. Remember that for $x \in \mathbb{R}$, there are \textit{no values} of $x$ that result in the exponential function $f(x)=\exp{x}$ being a non-positive real number.

Hankel \cite{Hankel}, Riemann \cite{riemann1859number}, and Edwards \cite{Edwards} all fail to provide any mathematically valid reason for equating the "strictly speaking" interpretation of the "first contour" on the left side of Eq. \ref{Hankel2} (the branch cut) to the "non-strictly speaking" interpretation of the "first contour" on the right side of Eq. \ref{Hankel2} (the Hankel contour). Again, the points on the contour represented by the left side of the equation (the branch cut) have no defined value, and thus are also non-holomorphic. As for the points on the Hankel contour represented by the right side of the equation, ("slightly above the Real axis as it descends from $+\infty$ to $0$ and slightly below the Real axis as it goes from $0$ back to $+\infty$"), they have defined values.
\footnote{But note: How far away from the branch cut do these points need to be in order to have defined values? Here we encounter the ancient "Sorites paradox", a.k.a "the paradox of the heap".}

Fortunately, in contrast to Riemann \cite{riemann1859number} and Edwards \cite{Edwards}, Whittaker \cite{Whittaker} \textit{does} provide a basis for equating the "strictly speaking" interpretation of the "first contour" on the left side of Eq. \ref{Hankel2} to the "non-strictly speaking" interpretation of the "first contour" on the right side of Eq. \ref{Hankel2}: \textit{the path equivalence corollary of Cauchy's integral theorem} is given as the mathematical basis for equating the Hankel contour to the branch cut. \footnote{See Whittaker \cite{Whittaker}, pp.85-7, 244, \S 5.2, Cor 1.}
However, this basis is \textit{neither mathematically nor logically valid}. The Hankel contour and the branch cut contradict the prerequisites of the Cauchy integral theorem, 
\footnote{See Whittaker \cite{Whittaker}, p.85.}  
and of its corollary. 
\footnote{See Whittaker \cite{Whittaker}, p.87.}  
Due to the LNC, these contradictions in the derivation of Riemann's Zeta function render it a paradox in "classical" logic, intuitionistic logic, and Zermelo-Fraenkel set theory, and therefore mathematically invalid.

This section presents the reasons why the Hankel contour contradicts the prerequisites of Cauchy's integral theorem. Cauchy's integral theorem states that if function $f(z)$ of complex variable $z$ is "holomorphic" (complex differentiable) at all points \textit{on} a simple closed curve ("contour") $C$, and if $f(z)$ is also holomorphic at all points \textit{inside} the contour, then the contour integral of $f(z)$ is equal to zero: 
\footnote{See Whittaker \cite{Whittaker}, p.85.}  
\begin{equation} \label{eq:2.1}
\int_{(C)} f(z)\cdot dz = 0
\end{equation}
The path equivalence corollary of Cauchy's integral theorem states the following: 
\footnote{See Whittaker \cite{Whittaker}, p.87, Cor 1.}  

(1) If there exist four distinct points ($z_0$, $Z$, $A$, and $B$) on the Cartesian plane (that represents the complex domain), and the two points $z_0$ and $Z$ are connected by two distinct paths $z_0AZ$ and $z_0BZ$ (one path going through $A$, the other path going through $B$), and 

(2) if function $f(z)$ of complex variable $z$ is holomorphic at all points on these two distinct paths $z_0AZ$ and $z_0BZ$, and $f(z)$ is holomorphic at all points enclosed by these two paths,  

(3) then any line integral connecting the two points $z_0$ and $Z$ inside this region (bounded by $z_0AZ$ and $z_0BZ$) has the same value, regardless of whether the path of integration is $z_0AZ$, or $z_0BZ$, or any other path disposed between $z_0AZ$ and $z_0BZ$. 

Riemann invalidly used Cauchy's integral theorem to assign, to the branch cut, the value of the limit of the Hankel contour (as the Hankel contour approaches the branch cut of $f(x)=\log(-x)$ at $x \in \mathbb{C}$). 

But by definition, $\log(-x)$ has no value (and thus is non-holomorphic) at all points on half-axis $x\in \mathbb{R}, x\ge0$. The geometric proof that $\log(-s)$ is non-holomorphic at all points on half-axis $s\ge0$ is as follows: In the Cartesian plane, the first derivative of $f(x)=\log(-x)$, for $x \in \mathbb{R}$ at a value of $x$, is represented by the slope of the line tangent to $f(x)$ at $x$. However, $f(x)$ has no values at $x\ge0$, so its first derivative cannot have any values at $x\ge0$. 

(Note however, that for $s \in \mathbb{C}$, there exists a definition for the branch cut of $f(s)=\log(-s)$
 that assigns to it the values of $f(s)=\log(|s|)$ (and remains undefined at $s=0$). This definition contradicts the definition of logarithms of Real numbers. 
\footnote{See the Encyclopedia of Mathematics \cite{EoM}.} 

Moreover, the Hankel contour is either open, or closed, at $x = +\infty$ (the latter enclosing non-holomorphic points). In both cases, the Hankel contour contradicts prerequisites of Cauchy's integral theorem.

If the Hankel contour is open, the Cauchy integral theorem cannot be used, because it only applies to closed contours. In the alternative, if the Hankel contour is indeed closed at $+\infty$ on the branch cut, as assumed by Riemann, 
\footnote{See Whittaker \cite{Whittaker}, p.245.} 
then the Hankel contour still contradicts the requirements of the Cauchy integral theorem. This is because the closed Hankel contour encloses the entire branch cut of $f(z)$, and the branch cut consists entirely of non-holomorphic points. Also, there would be a non-holomorphic point on the Hankel contour itself, at the point where it intersects the branch cut at $+\infty$ on the Real axis. These reasons disqualify the use of the Cauchy integral theorem with the Hankel contour.

For these reasons, it is not valid to use the Cauchy integral theorem's path equivalence corollary to find the limit of the Hankel contour, as the Hankel contour approaches the branch cut of $f(x)=\log(-x)$ at $x \in \mathbb{C}$. So the derivation of Riemann's Zeta function violates the Law of Non-Contradiction (LNC). 

\subsection{If Riemann's Zeta Function Were True, its Contradiction of Zeta's Dirichlet Series Would Create a Paradox}

Given that the Dirichlet series definition of the Zeta function is \textit{proven} to be divergent throughout the half-plane $\text{Re}(s)\le1$, if Riemann's Zeta function were true, then the Zeta function would have both a true  \textit{convergent} definition and a true \textit{divergent} definition throughout half-plane $\text{Re}(s)\le1$ (except at $s=1$). 

Moreover, if Riemann's Zeta function were true, and thus convergent throughout half-plane $\text{Re}(s)\le1$, then it would be \textit{convergent} throughout the Real half-axis $\{\text{Re}(s)<1, \text{Im}(s)=0\}$, which is a sub-set of the half-plane $\text{Re}(s)\le1$. (Riemann's functional equation of the Zeta function even claims to have "trivial zeros" on this Real half-axis.)
This result of "convergence" would directly contradict the results of "divergence" produced by the Integral test for convergence (a.k.a. the Maclaurin-Cauchy test for convergence) when applied to the Dirichlet series definition of the Zeta function, for all values of $s$ on this Real half-axis.

Also, if Riemann's Zeta function were true, and thus convergent throughout half-plane $\text{Re}(s)\le1$, then it would be \textit{convergent} at all points on the misleadingly-named "line of convergence" at $\text{Re}(s)=1$ (except at $s=1$). This would directly contradicts the divergence of the Dirichlet series definition of the Zeta function along this line. 
\footnote{See Hardy \cite{Hardy2}, p.5, Example (iii), citing Bromwich \cite{Bromwich}.} 

Each of these results would render the Zeta function a paradox in the half-plane $\text{Re}(s)\le1$, due to the contradictions. This, in turn, would be sufficient to cause "deductive explosion" for all other conjectures or theorems that would assume that the Zeta function were true in that half-plane.

However, according to the mathematical definitions of "convergence" and "divergence", a function cannot be \textit{both} convergent \textit{and} divergent at any value in its domain. 
\footnote{See Hardy \cite{Hardy}, p.1.}  
Moreover, if Riemann's definition were true, then the two contradictory definitions of the Zeta function in the half-plane $\text{Re}(s)\le1$ would also contradict the definition of a "function" in set theory (both naive and ZF), due to the one-to-two mapping from domain to range \footnote{See Stover \cite{Stover}.}  
Perhaps most alarmingly, if the two contradictory definitions of the Zeta function were both true in the half-plane $\text{Re}(s)\le1$, it would mean that the axiomatic system called "mathematics" would be inconsistent, thereby invalidating it according to logics having the Law of Non-Contraction (LNC) and the "Principle of Explosion" (ECQ)..

In addition, if both of the two contradictory definitions of the Zeta function were true, then this would violate all three of Aristotle's three "Laws of Thought". The contradictory definitions of Zeta would not only violate Aristotle's Law of Non-Contradiction (LNC), but also his Law of Identity (LOI) (according to which each thing is identical with itself), and his Law of the Excluded Middle (LEM) (because it would mean that at certain values of $s$, the Zeta function is simultaneously \textit{both} divergent \textit{and} convergent).  In summary, if both the Dirichlet series definition and Riemann's definition of Zeta are true, this result would violate all of the LOI, LEM, and LNC. The violation of LNC would cause ECQ ("Explosion").

Moreover, if Riemann's Zeta function were true, its violation of the LNC would mean that the foundation logic of mathematics (and therefore also of physics) \textit{could not} be Zermelo-Fraenkel set theory. Zermelo-Fraenkel set theory inherently has LNC and ECQ as axioms, because it was created in order to avoid the paradoxes of Frege's naive set theory (in particular, Russell's paradox).

Instead, the foundation logic of mathematics (and therefore also of physics) would have to be a paradox-tolerant logic, such as a three-valued logic (e.g. Bochvar's 3VL 
\footnote{See Bochvar \cite{Bochvar}.}
or Priest's $LP$ 
\footnote{See Priest \cite{Priest5}, \cite{Priest7}, and Hazen \cite{Hazen}.}),  
or a "paraconsistent" logic that has the LNC as an axiom, but not ECQ. 
\footnote{See Priest \cite{Priest1}, \cite{Priest3}, and \cite{Priest10}.} 

\subsection{Riemann's Zeta Function is False, so it Renders Unsound All Arguments that Falsely Assume it is True}

There is an error in the derivation of Riemann's Zeta function, 
\footnote{This is discussed in detail in section \ref{derivation} of this paper.}
due to Hankel's contour contradicting the preconditions of Cauchy's integral theorem. 
\footnote{This is fortunate, because otherwise the contradictory versions of Zeta would mean that mathematics is inconsistent, and thus invalid in logics that have LNC and ECQ.}

Yet even this result is problematic, because if Riemann's Zeta function is \textit{false} at all values of $s$ in half-plane $\text{Re}(s)\le1$ (except at $s=1$), then all mathematics conjectures and theorems, and physics theories, that falsely assume that Riemann's Zeta function is true are rendered \textit{unsound} (and invalid) in Aristotelian, classical, and intuitionistic logics (and even in the paradox-tolerant 3VLs and paraconsistent logics). 

For example, the "Zeta Function Regularization" used in physics is rendered invalid, because it equates a true definition of the Zeta function to a false definition. Moreover, because Riemann's Zeta function is one example of the Dirichlet $L$-functions, the falsity of Riemann's Zeta function is the example that disproves the assumption that all $L$-functions are true. More specifically, the false assumption that the $L$-functions are true includes the false assumption that Riemann's definition of "analytic continuation" 
is valid. 
\footnote{See Bruin \cite{Bruin}, p.4.}  
In turn, the false assumption that all $L$-functions are true renders unsound several mathematical theorems (e.g. the Modularity theorem, Fermat's last theorem) that are presumed to be proven.

The falsity of Riemann's Zeta function also confirms that $\zeta(1)\ne0$. This resolves the Birch and Swinnerton-Dyer (BSD) Conjecture in favor of finiteness. 
\footnote{See Clay \cite{Clay}.}  
The Dirichlet series exclusively defines the Zeta function, so at $s=1$, it is the "harmonic series", which is proven to be divergent by the Integral test for convergence 
\footnote{See Guichard  \cite{Guichard1}, Thm 13.3.4.}  
Moreover, the Landau-Siegel zero 
\footnote{See Siegel \cite{Siegel2}; and Conrey \cite{Conrey}, p.351.}  
is non-existent, due to the invalidity of of $L$-functions in general (resulting from the invalidity of Riemann's "analytic continuation" of the Zeta function). 

The falsity of Riemann's Zeta function, and of $L$-functions, 
\footnote{Such that the Zeta function is exclusively defined by Dirichlet series.}
resolves the BSD conjecture and triggers a "domino effect" of logical unsoundness (due to false assumptions) through a chain of conjectures that are proven to be equivalent. For example:
\begin{itemize}
    \item The BSD conjecture "for elliptic curves over global fields of positive characteristic" is equivalent to the Tate conjecture "for elliptic surfaces over finite fields", \footnote{See  Totaro \cite{Totaro}, citing Ulmer \cite{Ulmer2}, pp.6,31-32; Totaro \cite{Totaro2}, p.578;  and also Milne \cite{Milne4}, p.3, Thm 1.4.}
    \item The Tate conjecture is equivalent to the Hodge conjecture  "for abelian varieties of $CM$-type" \footnote{See Gordon \cite{Gordon}, p.364, \S 11.2, citing Pohlmann \cite{Pohlmann},  Piatetskii-Shapiro \cite{Piatetskii-Shapiro}, Borovoi \cite{Borovoi} and \cite{Borovoi2}; Deligne \cite{Deligne}, p.43, Cor 6.2.}  
    \footnote{See also Shioda \cite{Shioda}, p.60, citing Pohlmann \cite{Pohlmann}, \S 2,  Mumford \cite{Mumford}, Kubota \cite{Kubota},  Ribet \cite{Ribet},  Hazama \cite{Hazama}.}  
    \footnote{See also Beauville \cite{Beauville}, pp.12-14, Cor 5.5, citing Mattuck \cite{Mattuck}, Tate2 \cite{Tate2}, and Tankeev \cite{Tankeev}.}  
\end{itemize}

Therefore, the Tate conjecture and Hodge conjecture, which falsely assume that all $L$-functions are true, are rendered unsound by the falsity of Riemann's Zeta function in half-plane $\text{Re}(s)\le1$, via their relationships with the BSD conjecture.

There exist other conjectures rendered unsound by the falsity of Riemann's Zeta function in half-plane $\text{Re}(s)\le1$, due to their relationship to the BSD conjecture. These include the finiteness of the Tate-–Shafarevich group, and the finiteness of the Brauer group.  
\footnote{See Totaro \cite{Totaro2}, p.579; and Wiles \cite{Wiles}, p.2, citing: Tate \cite{Tate}, p.416,426; Milne \cite{Milne}, Cor 9.7.} 

Regarding Hadamard and de la Vallée Poussin’s respective proofs of the prime number theorem, Borwein (\cite{Borwein}) argues that they "follow from the truth of the Riemann hypothesis". 
\footnote{See Borwein \cite{Borwein}, pp.9,61, \S7.1, \S12.4; Edwards \cite{Edwards}, pp.68-69.} 
But Borwein is incorrect. Instead, these proofs are true because the Zeta function is exclusively defined by the Dirichlet series (which has no zeros). Therefore, the resulting Zeta function has no zeros on the misleadingly-named "line of convergence", $\text{Re}(s)=1$.

\subsection{A Third Definition of the Zeta Function}

\subsubsection{Derivation of the Third Definition}

Ash and Gross \cite{Ash} derive a third definition of the Zeta function from the original Dirichlet series definition. 
\footnote{See Ash and Gross \cite{Ash}, pp.169-171.}  
This third definition contradicts both Dirichlet’s and Riemann’s definitions of the Zeta function. 
Ash and Gross derive this definition of Zeta by multiplying the Dirichlet series of $\zeta(s)$ by the term $2^{-s}$:

\begin{equation}
\frac{1}{2^{s}}\cdot \zeta(s) = \frac{1}{2^{s}} + \frac{1}{4^{s}} + \frac{1}{6^{s}} + \frac{1}{8^{s}} + \cdots.
\end{equation}
(Note that this cannot be division by zero, because there is no value of $s$ for which $2^{(-s)}$ equals zero). This series is then twice subtracted from the original Dirichlet series, resulting in a conditionally convergent series:
\begin{equation}
(1 - \frac{1}{2^{s}} - \frac{1}{2^{s}})\cdot \zeta(s) =  1  - \frac{1}{2^{s}}  + \frac{1}{3^{s}}  - \frac{1}{4^{s}}  + \frac{1}{5^{s}}  - \frac{1}{6^{s}} + \cdots
\end{equation}
Note that the right side of the above equation is the Dirichlet series $\sum a_{n}n^{-s}$, wherein $a_n = (-1)^{n-1}$, so $|a_{1} + \cdots + a_{n}| < 2$ for all $n$. Rearranging the terms of the equation immediately above produces:
\begin{equation} \label{Ash_Zeta}
\zeta(s) = \Big(1-\frac{1}{2^{s-1}}\Big)^{-1}\cdot \Big(1 - \frac{1}{2^{s}} + \frac{1}{3^{s}} - \frac{1}{4^{s}} + \frac{1}{5^{s}} -
\frac{1}{6^{s}} + \cdots\Big)
\end{equation}
Ash and Gross \cite{Ash} cite the following theorem:
\footnote{See Ash and Gross \cite{Ash}, p.169, Theorem 11.7. The proof of this theorem can be found at Conrad \cite{Conrad} pp.2-3, Theorem 9, which cites Jensen \cite{Jensen} and Cahen \cite{Cahen}. See also Hardy \cite{Hardy2}, pp.3-5.} 
\begin{quotation}
Suppose that there is some constant $K$ so that $|a_{1} + \cdots + a_{n}| < K$ for all $n$. Then the Dirichlet series $\sum a_{n}n^{-s}$ converges if $\sigma > 0$.
\end{quotation}
and note that  $(1 - 1/2^{s} + 1/3^{s} - 1/4^{s} + 1/5^{s} -
1/6^{s} + \cdots)$ is a Dirichlet series that has the coefficients $1, -1, 1, -1, 1, \cdots$. So $|a_{1} + \cdots + a_{n}| < K$, and therefore the Dirichlet series $\sum a_{n}n^{-s}$ converges for $\sigma > 0$, as per Theorem 11.7.
This result proves that this third definition of the Zeta function is convergent throughout half-plane $\text{Re}(s)>0$ (except at $s=1$), and is divergent at the pole $s=1$ and throughout half-plane $\text{Re}(s)\le0$

However, note that the third definition of the Zeta function, as defined in Eq. \ref{Ash_Zeta} is \textit{absolutely} convergent for all values of $s$ in half-plane $\text{Re}(s)>1$ (which is consistent with the first two definitions of Zeta), but is \textit{conditionally} convergent throughout the "critical strip", $0<\text{Re}(s)\le1$. Therefore, according to the Riemann series theorem, it is a paradox in the critical strip.

\subsubsection{The Third Definition Contradicts the First Two Definitions}

This third definition of the Zeta function contradicts the Dirichlet series definition of the Zeta function throughout the “critical strip” ($0<\text{Re}(s)\le1$, except at the pole $s=1$), where the third definition is convergent and the Dirichlet series definition is divergent. The third definition contradicts Riemann’s Zeta function throughout half-plane ($\text{Re}(s)\le0$), where the third definition is divergent and Riemann's definition is convergent.

Clearly, in every logic that has any of Aristotle's three "Laws of Thought" as axioms (the Law of Identity (LOI), the Law of Non-Contradiction (LNC), and the Law of the Excluded Middle (LEM)), only \textit{one} of these three contradictory definitions of the Zeta function can be true. In such a logic, it is impossible for two or three of the definitions to be true. 

\subsubsection{In the Critical Strip, the Third Definition is a Conditionally Convergent Series (and is a Paradox There)}

The third definition of Zeta is \textit{conditionally} convergent in the "critical strip" ($0<\text{Re}(s)\le1$). Therefore, Riemann's series theorem proves that it can be rearranged to be divergent at domain values where it is \textit{conditionally} convergent. So the third definition of the Zeta function is both convergent and divergent in the "critical strip". This is a contradiction, and a paradox. 
\footnote{Series are classified into three categories (absolutely convergent, conditionally convergent, and divergent), so in a 3-valued logic, conditionally convergent series are assigned the third truth-value, because they are unlike the other two categories. }

In mathematics, this paradoxical result violates the definition of a "function" (due to the one-to-many mapping from domain to range), and also the associative and commutative properties of addition.
\footnote{Note that the two examples of conditionally convergent series shown in Weisstein \cite{Weisstein} are both represented as summations (of positive and negative numbers). These summations clearly show that conditionally convergent series contradict the associative and commutative properties of addition.}
In classical logic, this paradoxical result  violates the LOI, LNC, and LEM. The violation of LNC causes ECQ ("explosion"). In certain three-valued logics (e.g. Bochvar's 3VL), this paradox would be assigned the third truth-value ("paradox"), causing the LOI, LNC, and LEM to fail.

Moreover, all other \textit{conditionally} convergent series used in mathematics and physics theories violate the LNC, and cause ECQ (according to classical and intuitionistic logics). In certain three-valued logics (e.g. Bochvar's 3VL), these math and physics theories are assigned the third truth-value (e.g. "paradox" in Bochvar's 3VL).  

\subsubsection{The Zeta Function Has No Zeros, So the Riemann Hypothesis is a Paradox}

Both Riemann's definition and the third definition of the Zeta function are false. So the Zeta function is exclusively defined by the Dirichlet series definition. Therefore, the Zeta function has no zeros. This renders the Riemann hypothesis a paradox, due to "vacuous subjects". The Riemann hypothesis pertains to zeros that do not exist.

According to the logical concept of "material implication", and according to the logical theorem \textit{Ex Contradictione (Sequitur) Quodlibet} (ECQ, or "Explosion"), a false statement implies any other statement, true or false. In other words, in the proposition "If X then Y", if "X" is false, then regardless of whether "Y" is true or false, the proposition is true. In a more specific example, the Riemann hypothesis can be phrased as: "If the Zeta function equals zero, then its domain value is on the critical line." But we have shown that the Zeta function never equals zero. So in this example, X is false. Therefore, regardless of whether Y ("its domain value is on the critical line") is true or false, the proposition is true. So the Riemann hypothesis is a paradox. 

\section{Conclusion}

Quantum Electrodynamics (QED) renormalizaion is a paradox in "classical" logic, intuitionistic logic, and Zermelo-Fraenkel set theory. The Euler-Mascheroni constant is not a constant - it is a conditionally convergent series, and therefore is a paradox (according to the Riemann series theorem). Both QED the Euler-Mascheroni constant introduce contradictions into mathematical proofs, and therefore are mathematically invalid.  
Riemann's definition of the Zeta function is false, because its derivation relies upon the Hankel contour and Cauchy's integral theorem, but the Hankel contour contradicts the preconditions of Cauchy's integral theorem. So Zeta function regularization is logically and mathematically invalid.

Also, a third definition of the Zeta function, which contradicts both the Dirichlet series definition and Riemann's definition, is false. The Zeta function has no zeros, so the Riemann hypothesis is a paradox due to "vacuous subjects". 

Any mathematical or physical conjecture or theorem that assumes that any of the above-listed paradoxes or falsities are true is rendered invalid by contradiction, \textit{unless} the foundation logic of mathematics (and therefore also of physics) is a paradox-tolerant logic that rejects the Law of Non-Contradiction (LNC), the Principle of Explosion (ECQ), or both. In this scenario, the foundation logic cannot be the "classical logic" of Whitehead and Russell's \textit{Principia Mathematica}, or Heyting's formalization of Brouwer's intuitionistic logic, or Zermelo-Fraenkel (ZF) set theory, which inherently has the LNC and ECQ as axioms.  \footnote{ZF was created to avoid the paradoxes of Frege's naive set theory.}

\singlespacing
\bibliographystyle{acm}
\bibliography{sample.bib}

\end{document}